\documentclass[12pt,pra,aps,amssymb,amsfonts,amsmath,tightenlines]{revtex4}

\usepackage[dvips]{graphicx}
\usepackage{dsfont, color,txfonts} 

\newcommand{\re}{\mbox{$\rm e$}}
\newcommand{\rd}{\mbox{$\rm d$}}
\newcommand{\half}{\mbox{$\textstyle \frac{1}{2}$}}

\begin{document}

\title{Three candidate election strategy}
\author{Dorje C. Brody$^1$ and Tomooki Yuasa$^{2}$}

\affiliation{$^1$
School of Mathematics and Physics, University of Surrey, 
Guildford GU2 7XH, UK \\ 
$^2$ Faculty of Economics and Business Administration, Graduate School 
of Management, Tokyo Metropolitan University, Tokyo, 100-0005, Japan
}

\date{\today}

\begin{abstract}
The probability of a given candidate winning a future election is worked out in closed form as a function of (i) the current support rates for each candidate, (ii) the relative positioning of the candidates within the political spectrum, (iii) the time left to the election, and (iv) the rate at which noisy information is revealed to the electorate from now to the election day, when there are three or more candidates. It is shown, in particular, that the optimal strategy for controlling information can be intricate and nontrivial, in contrast to a two-candidate race. A surprising finding is that for a candidate taking the centre ground in an electoral competition among a polarised electorate, certain strategies are fatal in that the resulting winning probability for that candidate vanishes identically. 
\end{abstract}

\maketitle

%%%%%%%%%%%%

\section*{Introduction}

This paper is concerned with informational strategies behind an election 
involving three or more candidates. Suppose that there is an election coming 
up next year, and that the three candidates have the support rates of, say, 
38\%, 26\%, and 36\%, according to the current opinion poll statistics. Do 
these numbers represent the current likelihoods of the candidates winning 
the future election? If not, what are they? Can we derive a formula for a 
given candidate winning the election in the future? Because voter preferences 
change over time in accordance with information revealed to them, such 
a formula 
ought to be dependant on how information is managed between today and 
the future election day. In this paper, an exact formula will be worked out for 
the probability of winning a future election that depends explicitly on the 
model for the flow of information. 

With such a formula at hand, we are able to ask a range of strategic questions 
and find quantitative answers. For example, for the candidate lagging behind 
with only 26\% support today, in which way should they reveal policy 
information so as to maximise the realised probability of winning the future 
election? How would that differ from the strategy for other candidates? Are 
there advantages in the positioning of the candidate's political party within the 
political spectrum, for example, leaning further to the right or to the left? 
Our purpose here is to 
provide a systematic framework to answer questions of this kind, building 
on the previous work on two-candidate election models 
\cite{Brody1,Brody2,Brody3}. In particular, we shall reveal some 
surprising features that emerge from having more than two candidates in 
an electoral competition. 

It is worth remarking that there is a long history of building mathematical 
models to analyse various aspects of electoral competitions (see, for 
example, \cite{xxx,Shepsle,JH,MO,Coughlin,FP,FM,Lloyd,yyy,YAKM}). The 
models hitherto considered in the literature tend to be either deterministic, 
or else probabilistic but static. While the importance of the 
role played by information in electoral competitions is widely acknowledged, 
the models that have been proposed in the literature omit 
the impact of noise such as rumours, speculations, disinformation, and so 
on. In contrast, our formulation, building on a successful application of the 
noisy information-based approach to modelling the dynamics of financial 
assets \cite{BHM}, 
takes on board the revelation of noisy information, which in turn can be 
used to deduce the statistics of the outcomes of future elections. In this 
way, a sensitivity analysis can be pursued so as to determine the impact 
of different informational strategies on the future statistics. 

With these preliminaries, the present paper will be organised as follows. We 
begin with a brief introduction to the information-based approach to election 
modelling introduced in \cite{Brody1}. We then proceed by reviewing some of 
the key findings of \cite{Brody2} in which the probability of winning a future 
election is worked out in the case of a two-candidate electoral competition (or 
a ``yes-no'' referendum) for benefits of readers less acquainted with the 
material. We shall also indicate here how a large number 
of information sources can be aggregated in the form of a single information 
process, and how individual information source affects the overall information.
We then examine the case of an election involving three candidates, 
and work out explicit formulae for the probabilities of six different outcomes 
(orderings of candidates) in a future election. Our focus here will be a 
first-past-the-post electoral system, although the fact that we have explicit 
formulae for all different scenarios means that the approach can be used in 
other electoral systems such as a proportional representation system. We 
then explore the impact of candidate's positioning within the political spectrum. 
In this context it will be shown that among a polarised electorate the candidate 
with a centre-ground position can find that the probability of winning a future 
election is identically zero in a first-past-the-post electoral system, unless a 
significant amount of reliable information is revealed to create higher volatilities. 
We then discuss briefly how our model can be implemented 
in practice in an electoral competition. 
We conclude with a discussion on how closed-form formulae for the 
probabilities of the various different outcomes of an election can be obtained 
in the present scheme when there are more than three candidates.

\section*{Information-based modelling of electoral competition}

Let us begin by examining the role of information in modelling an electoral 
competition. In a generic situation, voters will have a range of issues that 
concern them in deciding which candidate to choose. The policy positions 
of the candidates, if elected,  on these issues, however, are known only 
partially to the voters. Nevertheless, as time progresses, more information 
about the candidates, or their views on different policy positions, are revealed, 
which in turn will shift voter perceptions. We can model this dynamics by use 
of the mathematical techniques of signal processing. Specifically, we model 
the position of candidate $l$ on issue $k$ by a random variable $X^{k,l}$, 
taking a range of values labelling different policy positions on this issue, 
where the probability that $X^{k,l}$ taking a given value represents the 
voter perception of the policy position of the candidate $l$ on this issue. 
Voters will then acquire partial information about the values of these 
random variables. For a fixed candidate and fixed issue, because the 
quantity of interest to voters is the candidate's position on that issue, 
represented by the value of $X^{k,l}$, we can think of this random variable 
as the `signal' in a communication channel, which is superimposed with 
noise that represents, for instance, rumours, speculations, disinformation, 
and so on. 
The arrival of information to voters can therefore be modelled in the form 
of a superposition of signal and noise. The idea we evoke here is that we 
regard the environment in which voters are living, in itself, as forming a 
communication channel, where noisy information is transmitted through 
newspapers, radio and television broadcasts, internet, word of mouth, and 
so on. Then based on this partial information, voters will come up with 
best estimates for each of the factors. Candidates 
are then ranked in accordance with these estimates, reflecting the 
preferences of voters. 

This concept has been formalised mathematically in \cite{Brody1} as a 
structural approach to modelling electoral competition, with an emphasis 
to determine the statistics of the impact of disinformation on democratic 
processes. For a better understanding of the 
information-based formalism it will be useful to explain the structural 
approach in more detail. The fundamental idea is to first model the 
information-providing process $\{\xi_t^{k,l}\}$ associated to the policy 
position $k$ of candidate $l$. For example, if the noise that obscures 
the value of $X^{k,l}$ is modelled by an additive Gaussian noise 
$\{B_t^{k,l}\}$, and if the rate at which information is revealed to the 
electorates at time $t$ is given by $\sigma_t^{k,l}$, then the information 
process takes the form of a `signal-plus-noise' decomposition: 
\begin{eqnarray}
\xi_t^{k,l} = X^{k,l} \int_0^t \sigma_s^{k,l} \, \rd s + B_t^{k,l}\, . 
\label{eq:1}
\end{eqnarray} 
If we let ${\cal F}_t$ denote the totality of information available to the 
electorates at time $t$ generated by (\ref{eq:1}) for all $k,l$, 
then the best estimate people will arrive at 
about the $k$-th policy position of candidate $l$ is given by the 
conditional expectation ${\hat X}_t^{k,l}={\mathbb E}[X^{k,l}|
{\cal F}_t]$. 

Next, for a given voter $m$, we let $w_{m,k}$ denote the preference 
weight of that voter for issue $k$. Then we can identify the `score' 
$S_m^l(t)$ assigned at time $t$ by voter $m$ for candidate $l$. For 
example, in a linear scoring system we have 
\begin{eqnarray}
S_m^l(t) = \sum_k w_{m,k} \, {\hat X}_t^{k,l}. 
\end{eqnarray} 
Then at time $T$ of the election, voter $m$ will choose the candidate 
with the highest score at that time. Because a large number of the voter 
preferences $\{w_{m,k}\}$ can be sampled from a distribution 
\cite{Brody1}, in a structural approach it is possible to engage in a rather 
detailed issue-by-issue scenario analysis to identify optimal informational 
strategies, as well as making predictions of the statistics of a future 
election. 

\section*{Reduced-form approach to electoral competition}

An alternative `reduced form' approach has also been introduced in 
\cite{Brody1} and was further developed in \cite{Brody2}. 
The idea of a reduced form approach is to aggregate a broad range of 
issues into a single random variable $X$ that labels different candidates. 
We remark that the terminologies of structural versus 
reduced-form are derived from an analogous consideration in the 
modelling of credit risk in financial markets. Here, for a given cash 
flow associated to a given firm, one can attempt to go into a detailed 
structural analysis of that firm in identifying the risk associated to that 
cash flow. However, in most cases this is not feasible because 
relevant financial 
structures -- for instance other cash flows linked to that firm -- are far 
too complicated and often not even transparent. To remedy this issue, 
in credit risk modelling an alternative reduced-form approach has been 
introduced as a way for capturing the essence of credit 
risks without going into any of the structural details so that the method 
can be applied universally and systematically, in practical applications. 

In contrast to credit risk modelling, for an electoral 
competition a structural approach is entirely feasible. This is because 
the number of significant issues that are of concern to a large number 
of electorates is limited -- perhaps a handful as opposed to thousands 
in the case of a typical credit product -- and likewise the number of 
candidates is also limited. Nevertheless,
the advantage of the reduced form formalism, which 
will be explained in more detail now, is that 
just as in credit risk modelling, the mathematical analysis 
can be simplified considerably without losing many of the key features of 
the structural formalism. In the present paper we shall therefore develop 
the theory underlying a three candidate race in the reduced-form approach. 

In a reduced form approach to modelling electoral competitions, 
a wide-ranging information relevant to deciding which candidate to choose 
is aggregated in the form of a single information process that represents 
abstractly the choice of the candidates. More specifically, in an election 
with $N$ candidates we let $X$ be the random variable taking the values 
$\{x_k\}$, $k=1,2,\ldots,N$, that label different candidates, with the 
probabilities $\{p_k\}$. These probabilities represent the current opinion 
poll statistics. These statistics will change in time, in accordance with the 
revelation of information related to the candidates. Recall that voters 
wish to determine which candidate to vote for. Hence in a reduced form 
model the random variable $X$ plays the role of the `signal'. 

As an elementary model let us assume that information concerning the 
candidates is revealed to the voters at a constant rate $\sigma$, and that 
wide-ranging noise is modelled by a Brownian motion $\{B_t\}$. Then 
the flow of information takes the familiar signal-plus-noise form 
\begin{eqnarray}
\xi_t=\sigma X t + B_t . 
\label{eq:x1} 
\end{eqnarray}
In a more realistic scenario the information flow rate is time dependent, 
and in this case the signal component is modified to $X\int_0^t\sigma_s 
\,{\rm d}s$. However, for simplicity we shall consider the special case in 
which $\sigma_t=\sigma$ is constant, with the remark that all the results 
presented below can be extended to the time-dependent case without 
losing analytical tractability. 

It is worth remarking that in a real election there are many 
information sources associated to the candidate choice $X$. We thus 
have the information process $\xi_t^j = \sigma_j X t + B_t^j$ 
generated by the $j$-th information source. However, without loss of 
generality we can aggregate different information sources in the form 
of a single information process (\ref{eq:x1}). To see this, let us first 
consider a pair of information sources: $\xi_t^1 = \sigma_1 X t + B_t^1$ 
and $\xi_t^2 = \sigma_2 X t + B_t^2$. Letting $\rho={\mathbb E}
[B_t^1 B_t^2]$ denote the correlation of the two Brownian noise, it is 
an elementary fact that there exists a Brownian motion ${\bar B}_t$, 
independent of $B_t^1$, such that we have $B_t^2 = \rho B_t^1 + 
\sqrt{1-\rho^2} {\bar B}_t$. Defining 
\begin{eqnarray} 
{\bar\sigma}=\frac{\sigma_2-\rho\sigma_1}{\sqrt{1-\rho^2}} 
\quad {\rm and} \quad {\bar\xi}_t = {\bar\sigma}Xt+{\bar B}_t 
\end{eqnarray}  
we deduce that the information generated jointly by $\xi_t^1$ and 
$\xi_t^2$ is equivalent to that generated jointly by $\xi_t^1$ and 
${\bar\xi}_t$. Next, let us define 
\begin{eqnarray}
\sigma^2 = \sigma_1^2+{\bar\sigma}^2 \quad {\rm and} \quad 
B_t = \frac{1}{\sigma} (\sigma_1 B_t^1 + {\bar\sigma} {\bar B}_t) , 
\end{eqnarray}
and set 
\begin{eqnarray}
\xi_t = \sigma X t + B_t  \quad {\rm and} \quad 
\delta_t = \frac{\xi_t^1}{\sigma_1} - \frac{{\bar\xi}_t}{{\bar\sigma}} = 
\frac{B_t^1}{\sigma_1} - \frac{{\bar B}_t}{{\bar\sigma}} .  
\end{eqnarray}
Then a short calculation shows that the information generated jointly 
by $\xi_t^1$ and $\xi_t^2$ is equivalent to that generated jointly by 
$\xi_t$ and $\delta_t$, and that $\xi_t$ and $\delta_t$ are independent. 
However, $\delta_t$ is independent of the 
random variable $X$. It follows that ${\mathbb E}[X|\{\xi_t^1,\xi_t^2\}] 
= {\mathbb E}[X|\{\xi_t\}]$. In other words, the aggregate of the two 
information sources can be represented in the form of a single 
information process (\ref{eq:x1}). 

More generally, suppose that we have a series of information processes 
of the form 
\begin{eqnarray}
%(\eighthnote)\quad 
\left\{ \begin{array}{l} 
\xi^1_t = \sigma_1 X t + B^1_{t} \\ 
\quad  \,\,\, \vdots  \\ 
\xi^n_t = \sigma_n X t + B^n_{t} 
\end{array} \right. \label{eq:3}
\end{eqnarray}
available to the electorates about the choice of candidates, where the 
various noise processes $\{B_{t}^i\}_{i=1,\ldots,n}$ in general may be 
mutually correlated, with the correlation matrix $\rho_{ij}$. Then the 
aggregate of the $n$ information processes (\ref{eq:3}) can be 
represented by means of a single information process of the form 
(\ref{eq:x1}), where 
\begin{eqnarray}
\sigma^2 =\frac{\sum^n_i \sigma_i^2\, \rho^{-1}_{ii} -2 
\sum_{i \neq j}\sigma_i \, \sigma_j \, \rho_{ij}^{-1}}{\det(\rho)} 
\label{eq:4}
\end{eqnarray}
represents the \textit{effective information flow rate} and
\begin{eqnarray}
B_{t}=\frac{1}{\sigma}\left( \sum^n_{i,j} \sigma_i \, 
\rho^{-1}_{ij}\, B_{t}^i \right)
\end{eqnarray}
represents the \textit{effective noise}. 
Here $\rho^{-1}_{ij}$ denotes the $ij$ element of the 
inverse correlation matrix. This observation shows that although the 
idea of representing a wide range of information flows in terms of a 
single information process at first may seem restrictive, in fact it is 
quite general. 

Given our model (\ref{eq:x1}) to represent the flow of information, the 
initial voter preference for a candidate, embodied in the \textit{a priori} 
probability $p_i = {\mathbb P}(X=x_i)$, will change into the posterior 
preference $\pi_{it} = {\mathbb P}(X=x_i|\{\xi_s\}_{0\leq s\leq t})$, which 
is just the conditional probability that candidate labelled by $x_i$ should 
be chosen, given the information available up to time $t$. In the present 
example in which the information flow rate is constant, the information 
providing process $\{\xi_t\}$ is Markov, from which it follows that the 
conditional probability simplifies into $\pi_{it}={\mathbb P}(X=x_i|\xi_t)$. 
Then by use of the Bayes formula 
\begin{eqnarray}
{\mathbb P}(X=x_i|\xi_t) &=&  \frac{{\mathbb P}(X=x_i)\rho
(\xi_t|X=x_i)}{\sum_{j} {\mathbb P} (X=x_j)
\rho(\xi_t|X=x_j)} ,
\label{eq:x2}
\end{eqnarray} 
along with the fact that the conditional density function $\rho(\xi_t|X=x_i)$ 
for the random variable $\xi_t$ is Gaussian and is given by  
\begin{eqnarray}
\rho(\xi|X=x_i)=\frac{1}{\sqrt{2\pi t}} \exp\left(-
\frac{(\xi-\sigma x_i t)^2}{2t}\right),  
\label{eq:x3} 
\end{eqnarray} 
we deduce at once that 
\begin{eqnarray}
\pi_{it} =\frac{p_i\exp\left( \sigma x_i \xi_t-
\frac{1}{2} \sigma^2 x_i^2 t\right)} {\sum_j p_j 
\exp\left( \sigma x_j \xi_t-\frac{1}{2} \sigma^2 x_j^2 t\right)} .
\label{eq:x4}
\end{eqnarray}
Incidentally, this expression in the context of signal 
processing is known as the Wonham filter \cite{wonham} associated to 
the detection of a random drift of a Brownian motion; a problem that 
has also been explored more recently in different contexts 
\cite{Buonaguidi,Ekstrom}.

\section*{Two-candidate scenario}

Given the expression (\ref{eq:x4}) for the \textit{a posteriori} probability 
that the $i$th candidate should be chosen, we are able to ask a range 
of questions that link informational strategies to the election outcome. To 
this end we determine first the probability that the $i$th candidate wins 
a future election, given (a) current support rate, and (b) how information 
is managed from now to the election day, in the case of a two-candidate 
competition. In this case we may let, without loss of generality, the random 
variable $X$ labelling the two candidates be binary, taking the values $0$ 
and $1$. Let $p$ be the probability that $X=0$ and $1-p$ be the probability 
that $X=1$. Then the \textit{a posteriori} probability that, say, candidate $0$ 
being chosen when an arbitrary voter is sampled from the population is 
\begin{eqnarray}
\pi_{0t} = \frac{p}{p + (1-p) \exp\left( \sigma \xi_t-\frac{1}{2} \sigma^2 t\right)} . 
\label{eq:x5} 
\end{eqnarray} 
If the election were to take place in $T$ years time from today 
(today always implies time $t=0$), 
then the probability that candidate zero should win the election, subject to 
current poll and how information is revealed from today to the election day, 
is therefore given by ${\mathbb P}(\pi_{0T}>1/2)$. 

It is important to note that the realised winning probability ${\mathbb P}
(\pi_{0T}>1/2)$ for candidate zero, as of today, can be very different from 
the current support rate $p$. To understand this, imagine that the election 
is to take place in a week, and that candidate zero currently has $55\%$ 
support. Then unless something radical happens --- such as a revelation 
of a major scandal --- it is likely that voter preferences will not change 
very much in one week, hence candidate zero will receive nearly $55\%$ 
of the votes to secure a victory. In other words, in such a scenario the 
realised probability of candidate zero winning the future election is close 
to $100\%$, even though the support rate remains only $55\%$. Putting 
the matter differently, today's poll statistics is not the predictor for the 
likelihood of winning a future election, although it can be used to calculate 
the likelihood. 

One advantage of our approach is that we are able to derive an explicit 
formula for the probability of a given candidate winning the election 
that reflects this observation. The 
detailed derivation of the formula in the case of a two-candidate election 
is provided in \cite{Brody2}, which we shall not reproduce here. 
Instead, we mention some key steps that will be exploited in the present 
context. Namely, that the denominator of the conditional probability 
(\ref{eq:x4}) can be used to change probability measure ${\mathbb P}$ 
into a new measure ${\mathbb Q}$ such that under ${\mathbb Q}$ the 
information process $\{\xi_t\}$ is a standard Brownian motion \cite{Brody2}. 
Specifically, writing 
\begin{eqnarray}
\Phi_t =   p + (1-p) \exp\left( \sigma \xi_t-\half \sigma^2 t\right)  
\end{eqnarray}  
for the measure-change martingale, we have 
\begin{eqnarray}
{\mathbb P}\left(\pi_{0T}>\half\right) = {\mathbb E}^{\mathbb P}
\left[ {\mathds 1}\left\{ \pi_{0T}>\half\right\} \right] = {\mathbb E}^{\mathbb Q} 
\left[ \Phi_T \, {\mathds 1}\left\{ \pi_{0T}>\half \right\} \right] . 
\label{eq:x7}
\end{eqnarray}
Now the condition that $\pi_{0T}>1/2$ is equivalent to the condition that 
\begin{eqnarray}
\frac{\xi_T}{\sqrt{T}} < \frac{ \log\left( \frac{p}{1-p}\right) + 
\frac{1}{2}\sigma^2 T} {\sigma \sqrt{T}} , 
\label{eq:x8} 
\end{eqnarray} 
but under ${\mathbb Q}$ the information process is a Brownian motion, and 
hence $\xi_T/\sqrt{T}$ is a standard normal random variable. It then follows 
at once that 
\begin{eqnarray} 
{\mathbb P}\left(\pi_{0T}>\half\right) = p\, N(d^+) + (1-p)\, N(d^-) , 
\label{eq:x9}
\end{eqnarray} 
where 
\begin{eqnarray}
d^\pm =  \frac{ \log\left( \frac{p}{1-p}\right) \pm \frac{1}{2}\sigma^2 T}
{\sigma \sqrt{T}} 
\label{eq:x10} 
\end{eqnarray} 
and 
\begin{eqnarray} 
N(x) =  \frac{1}{\sqrt{2\pi}} \int_{-\infty}^{x}  \re^{-\frac{1}{2}z^2} \rd z  
\end{eqnarray} 
denotes the cumulative normal distribution function. 
It is a curious coincidence that the winning probability 
of a candidate in a two-candidate election is essentially the same as 
the option pricing formula of Black and Scholes in financial modelling. 
The formula shows, for instance, that if candidate zero has $55\%$ 
support rate today and if the election is to take place in a week, then 
even if the information revelation rate is as large as, say, $\sigma=1.2$, 
the winning probability will be about $89\%$; whereas if the information 
revelation rate is reduced to, say, $\sigma=0.5$, then this probability 
increases to $99.8\%$. In other words, the model reflects the intuition 
described above. Putting it differently, formula (\ref{eq:x9}) allows 
us to interpolate between today's and future's statistics.

\begin{figure}[t!]
\centerline{
\includegraphics[width=0.75\textwidth]{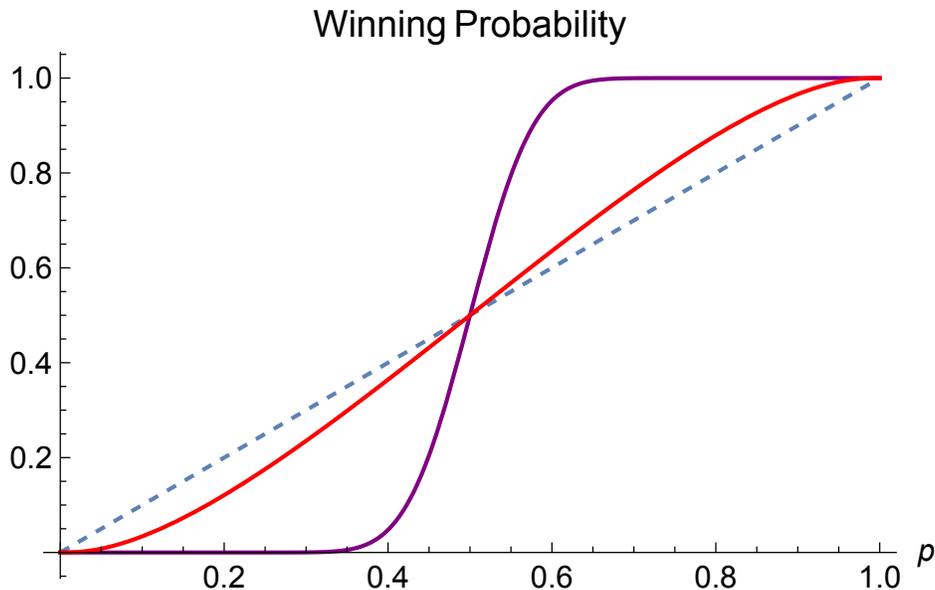}
}
\caption{\footnotesize{
\textit{Winning likelihood}. 
The probability that candidate zero will win the election in eighteen months 
($T=1.5$), as a function of the current support rate $p$ for the candidate. 
The realised likelihood of winning a future election is always higher than 
today's poll if $p>\frac{1}{2}$; and conversely lower than the poll if 
$p<\frac{1}{2}$. How much the winning probability deviates from the current 
poll depends on how much information is revealed over the next eighteen 
months. Here, two examples are shown, corresponding to the values 
$\sigma=0.2$ (in purple) and $\sigma=1.2$ (in red). 
}}
\label{fig1}
\end{figure}

We therefore see how an explicit formula (\ref{eq:x9}) for a given candidate 
winning a future election can be obtained in the case of a two-candidate 
competition. The winning probability, more explicitly, depends on the 
following three ingredients: (i) The current support rate $p$ for the candidate, 
(ii) the time $T$ left to the election, and (iii) the rate $\sigma$ at which 
information is revealed to the electorate from now until the election day. The 
only `control' parameter at candidates' disposal therefore is the information 
flow rate $\sigma$. To gain a better intuition here, therefore, let us examine 
how the winning probability depends on the current poll statistics for 
different values of $\sigma$. In Figure~\ref{fig1} we plot the winning 
probability for candidate zero as a function of the current support rate $p$ 
for two different values of $\sigma$. It is evident that if very little information 
is revealed from today to the election day, then the current support rate will 
not change significantly so that the likelihood of winning the election is 
considerably higher (lower) than the current poll if it is higher (lower) than 
$50\%$. What this means is that in a two-candidate election, if a candidate 
is losing then it is in their interest to release as much information as possible 
to generate more volatility; 
whereas if the candidate is winning then it is in their interest to conceal as 
much information as possible. This situation may be empirically familiar to 
some election strategists. When there are three or more candidates, 
however, there are some nontrivial situations that can arise, as we shall 
discuss below. 

\section*{Three-candidate electoral competition}

In a reduced-form approach, when there are three candidates we let the 
random variable $X$ take the values $\{x_k\}$ $(k=1,2,3)$, labelling the 
these candidates. Thus 
the event that $X=x_3$, say, means that the third candidate is being selected 
as the `right' choice, \textit{if the voters were to have access to the information 
$\{\xi_t\}$ eternally}. That is, in this scenario, which occurs 
with probability $p_3$, the sample path of $\pi_{3t}$ will be such that 
we have $\pi_{3t}\to1$ as 
$t\to\infty$. Of course, the election will take place earlier, so the vote share 
for the third candidate on the election day will instead be $\pi_{3T}$, which 
may or may not be larger than the support rates of the other two candidates.  
Hence even in the event in which $X=x_3$, in general this has little to do 
with the likelihood of the third candidate winning the election, unless the 
value of $\sigma$ is unusually large. 

As a matter of interpretation, to further clarify the meaning of the random 
variable $X$, we remark that our model setup is such that the value of 
$X$ will only be revealed base on the information process $\{\xi_t\}$ over 
an infinite time horizon, because 
\begin{eqnarray} 
\lim_{t\to\infty} \frac{\xi_t}{\sigma t} = X 
\end{eqnarray} 
in a distributional sense. In other words, if hypothetically the electorates 
were to live forever, then they will all learn which of the candidates they 
should all be voting for. However, the election will take place considerably 
sooner at time $T$, at which point there has not been sufficient information 
generated by the electorates according to the information process (unless 
$\sigma$ is very large). Thus the voter variability remains high, and the 
election outcome remains uncertain. It follows that our model is only of 
use until the election day, at which point it can be discarded, but this is 
all we need in order to interpolate the statistics between today and the 
election day. 

Differently stated, in a typical election cycle the information-providing 
process $\{\xi_t\}$ 
ceases to exist after the election day (or, equivalently, the information flow 
rate parameter $\sigma$ goes to zero, leaving behind nothing but noise), 
so the voters will never learn which candidate would have been the 
`right' candidate. Hence in reality none of the \textit{a posteriori} 
probabilities will 
actually converge to unity (except perhaps in certain extreme 
circumstances). This is because no one will have access to sufficient 
information to allow them to appropriately assign probabilities on 
counterfactual events: What would have happened if that person were 
elected rather than this person. Indeed, even in an event whereby a 
candidate who lost the election remains 
engaged towards a subsequent election, it is not always the case that much 
information about that candidate reaches the voters, as exemplified by an 
18 August 1996 Newsweek article titled ``Forgotten, but not gone'' about 
the then presidential candidate Ross Perot. 

In the case of a three-candidate race, there are six possible ordering scenarios 
for the vote share. Hence in order to determine the probability that the 
third candidate, say, wins the election, we need to work out the probabilities 
for individual ordering. This follows because the current probability for the 
third candidate to win the election is given by ${\mathbb P}(\pi_{3T}>\pi_{2T}>
\pi_{1T})+{\mathbb P}(\pi_{3T}>\pi_{1T}>\pi_{2T})$; and similarly for the 
other candidates to win the election. To this end we note from (\ref{eq:x4}) 
that for any $j\neq k$ the event $\pi_{kT}>\pi_{jT}$ holds true if and only if 
\begin{eqnarray} 
p_k\exp\left( \sigma x_k \,\xi_T-\half \sigma^2 x_k^2 T\right) >  
p_j\exp\left( \sigma x_j \,\xi_T-\half \sigma^2 x_j^2 T\right) . 
\label{eq:x12} 
\end{eqnarray} 
This condition can easily be solved for a condition on $\xi_T$, provided 
that the ordering of the numbers $\{x_i\}$ are given. Without loss of 
generality we may assume that $x_3>x_2>x_1$. Then from (\ref{eq:x12}) 
it follows that $\pi_{kT}>\pi_{jT}$ holds if and only if 
\begin{eqnarray} 
\xi_T > \frac{\log(p_j/p_k)+\frac{1}{2}\sigma^2(x_k^2-x_j^2)T}{\sigma(x_k-x_j)} ,
\label{eq:x13} 
\end{eqnarray} 
provided that $x_k>x_j$. Otherwise, the inequality is reversed. 

To proceed, let us introduce the notation 
\begin{eqnarray}
z_{kj} = \frac{\log(p_j/p_k)+\frac{1}{2} 
\sigma^2(x_k^2-x_j^2)T}{\sigma(x_k-x_j)} . 
\end{eqnarray} 
Then a short calculation shows, on account of the symmetry property 
$z_{kj}=z_{jk}$, that 
\begin{align*}
\pi_{3T}<\pi_{2T}<\pi_{1T} \quad &\Leftrightarrow \quad \xi_{T}
<\min\{z_{12},z_{23}\}, \\
\pi_{2T}<\pi_{3T}<\pi_{1T} \quad &\Leftrightarrow \quad 
z_{23}<\xi_{T}<z_{31} \qquad (\text{if}~ z_{23}<z_{31}), \\
\pi_{3T}<\pi_{1T}<\pi_{2T} \quad &\Leftrightarrow \quad 
z_{12}<\xi_{T}<z_{31} \qquad (\text{if}~ z_{12}<z_{31}), \\
\pi_{1T}<\pi_{3T}<\pi_{2T} \quad &\Leftrightarrow \quad 
z_{31}<\xi_{T}<z_{23} \qquad (\text{if}~ z_{31}<z_{23}), \\
\pi_{2T}<\pi_{1T}<\pi_{3T} \quad &\Leftrightarrow \quad 
z_{31}<\xi_{T}<z_{12} \qquad (\text{if}~ z_{31}<z_{12}), \\
\pi_{1T}<\pi_{2T}<\pi_{3T} \quad &\Leftrightarrow \quad \max\{z_{12},z_{23}\}<\xi_{T}.
\end{align*} 
Note here that, for example, the event $\pi_{2T}<\pi_{3T}<\pi_{1T}$ cannot 
be realised if $z_{23}>z_{31}$, and similarly for other three intermediate 
cases. With these conditions at hand, let us note that the probability 
${\mathbb P}(a<\xi_T<b)$ for any $a<b$ can be worked out by changing 
the probability measure. Specifically, we use the common denominator 
\begin{eqnarray}
\Phi_t = \sum_{j=1}^3 p_j 
\exp\left( \sigma x_j \xi_t-\half \sigma^2 x_j^2 t\right) 
\end{eqnarray} 
to effect a measure change ${\mathbb P}\to{\mathbb Q}$ so that under 
${\mathbb Q}$ the information process $\{\xi_t\}$ is a standard Brownian 
motion. Then we have 
\begin{eqnarray}
{\mathbb P}\left(a<\xi_{T}<b\right) &=& 
{\mathbb E}^{{\mathbb Q}}\left[\Phi_{T}{\mathds 1}\{a<\xi_{T}<b\}\right] 
\nonumber \\ &=& 
{\mathbb E}^{{\mathbb Q}}\left[\sum_{j=1}^{3}p_{j}\exp\left( 
\sigma x_{j}\xi_{T}-\half\sigma^{2}x_{j}^{2}T\right) 
{\mathds 1}\{a<\xi_{T}<b\}\right] \nonumber \\ &=& 
\sum_{j=1}^{3}p_{j}\left[ N\left(\frac{b-\sigma x_{j}T}{\sqrt[]{T}}\right) 
-N\left(\frac{a-\sigma x_{j}T}{\sqrt[]{T}}\right)\right] .
\end{eqnarray}

\begin{figure}[t!]
\centerline{
\includegraphics[width=0.3\textwidth]{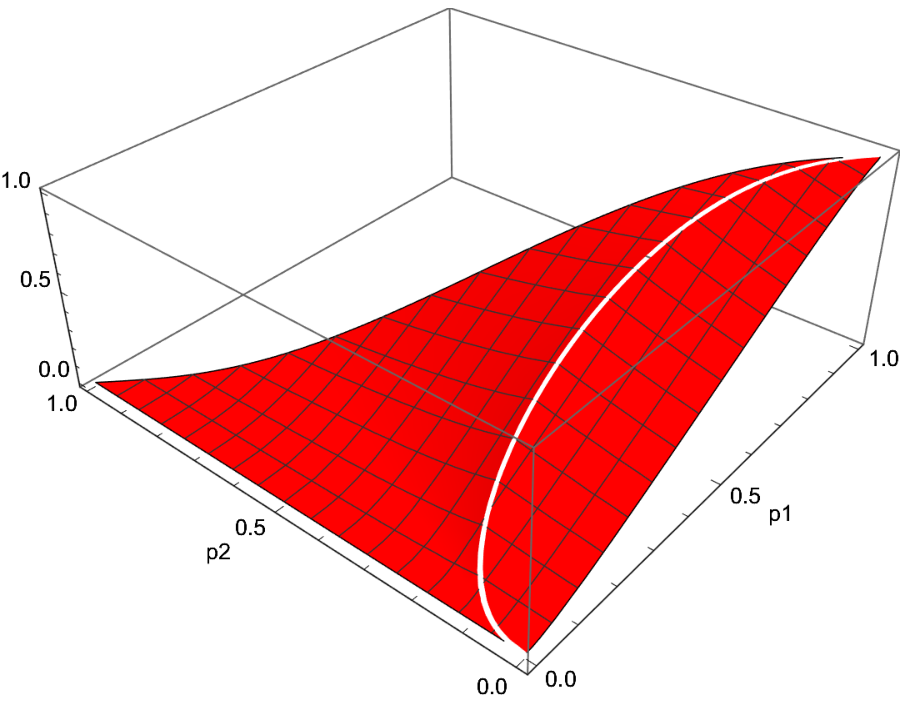}
\includegraphics[width=0.3\textwidth]{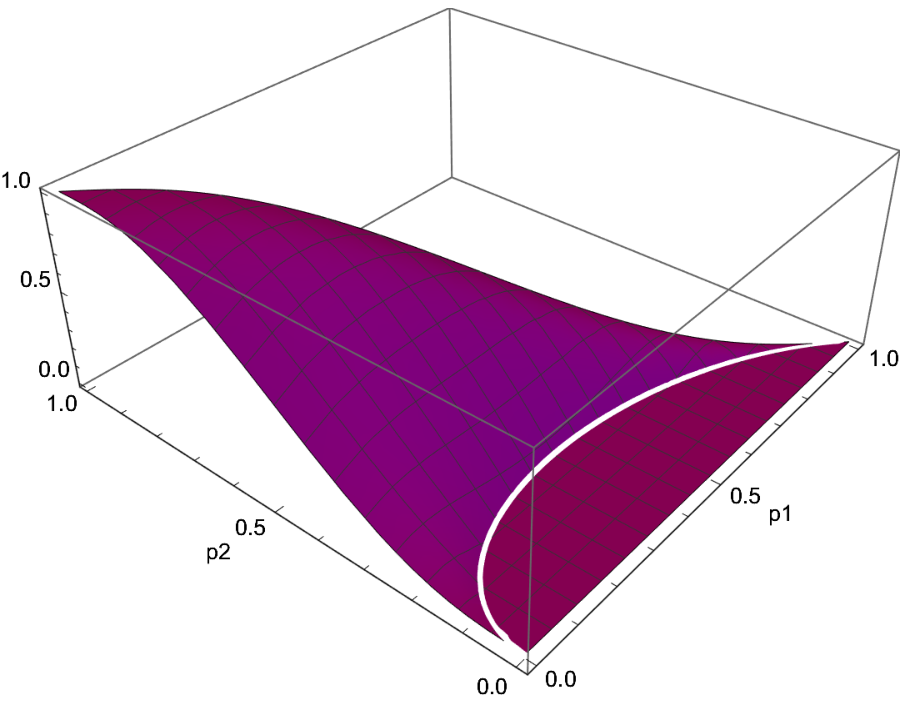}
\includegraphics[width=0.3\textwidth]{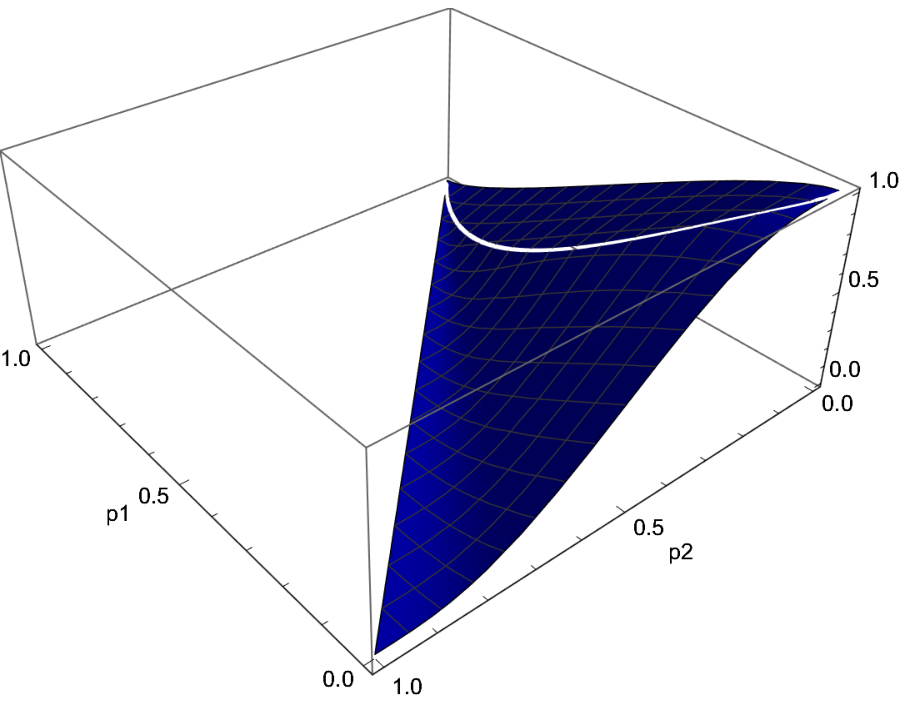}
}
\caption{\footnotesize{
\textit{Winning probabilities as functions of $(p_1,p_2)$}. The probabilities 
of winning a future election to take place in one year time ($T=1$), when 
the information flow rate is set at $\sigma=1$, are plotted here for the 
parameter choice $(x_1,x_2,x_3)=(1,2,3)$. The forms of the probabilities 
for candidate 1 (left panel, in red) and candidate 3 (right panel, in blue) are 
entirely symmetric. However, the behaviour of the probability for candidate 
2 (centre panel, in purple) 
is slightly different in that there is a region in the parameter space 
$(p_1,p_2)$ of the current support rates for which the probability of 
candidate 2 winning is identically zero. We will have more to say about this 
in the next section. 
}}
\label{fig:2}
\end{figure}

With these results at hand we are able to work out the probabilities for the 
six possible outcomes. For concreteness let us write them down explicitly 
here. They are 
\begin{eqnarray}
{\mathbb P}(\pi_{3T}<\pi_{2T}<\pi_{1T})
=\sum_{j=1}^{3}p_{j}\left[N\left(\frac{z_{12}-\sigma x_{j}T}{\sqrt[]{T}}\right)
{\mathds 1}\{z_{12}<z_{23}\}+N\left(\frac{z_{23}-\sigma x_{j}T}{\sqrt[]{T}}\right)
{\mathds 1}\{z_{12} \geq z_{23}\}\right] 
\end{eqnarray}
and
\begin{eqnarray}
{\mathbb P}(\pi_{2T}<\pi_{3T}<\pi_{1T})=\sum_{j=1}^{3}p_{j}\left[
N\left(\frac{z_{31}-\sigma x_{j}T}{\sqrt[]{T}}\right)-N\left(
\frac{z_{23}-\sigma x_{j}T}{\sqrt[]{T}}\right)\right]{\mathds 1}\{z_{23}<z_{31}\} ,
\end{eqnarray}
together determine the probability of the first candidate winning, 
\begin{eqnarray}
{\mathbb P}(\pi_{3T}<\pi_{1T}<\pi_{2T})=\sum_{j=1}^{3}p_{j}\left[ 
N\left(\frac{z_{31}-\sigma x_{j}T}{\sqrt[]{T}}\right)-
N\left(\frac{z_{12}-\sigma x_{j}T}{\sqrt[]{T}}\right)\right]
{\mathds 1}\{z_{12}<z_{31}\}
\end{eqnarray}
and 
\begin{eqnarray}
{\mathbb P}(\pi_{1T}<\pi_{3T}<\pi_{2T})=\sum_{j=1}^{3}p_{j}\left[ 
N\left(\frac{z_{23}-\sigma x_{j}T}{\sqrt[]{T}}\right)-
N\left(\frac{z_{31}-\sigma x_{j}T}{\sqrt[]{T}}\right)\right]
{\mathds 1}\{z_{31}<z_{23}\}, 
\end{eqnarray}
together determine the probability of the second candidate winning, and 
\begin{eqnarray}
{\mathbb P}(\pi_{2T}<\pi_{1T}<\pi_{3T})=\sum_{j=1}^{3}p_{j}\left[
N\left(\frac{z_{12}-\sigma x_{j}T}{\sqrt[]{T}}\right)-
N\left(\frac{z_{31}-\sigma x_{j}T}{\sqrt[]{T}}\right)\right]
{\mathds 1}\{z_{31}<z_{12}\}
\end{eqnarray}
and
\begin{eqnarray}
{\mathbb P}(\pi_{1T}<\pi_{2T}<\pi_{3T})=\sum_{j=1}^{3}p_{i}\left[
N\left(-\frac{z_{23}-\sigma x_{j}T}{\sqrt[]{T}}\right)
{\mathds 1}\{z_{12}<z_{23}\}+
N\left(-\frac{z_{12}-\sigma x_{j}T}{\sqrt[]{T}}\right)
{\mathds 1}\{z_{12} \geq z_{23}\}\right],
\end{eqnarray}
together determine the probability of the third candidate winning. In this way 
we obtain an explicit formula for each of the candidate winning the future 
election, as functions of (i) the current support rates $\{p_j\}$ for the 
candidates, (ii) the time $T$ left to the election, (iii) the rate $\sigma$ at 
which information is revealed to the electorate, and (iv) the choice of the 
candidate labels $\{x_j\}$. 

\begin{figure}[t!]
\centerline{
\includegraphics[width=0.45\textwidth]{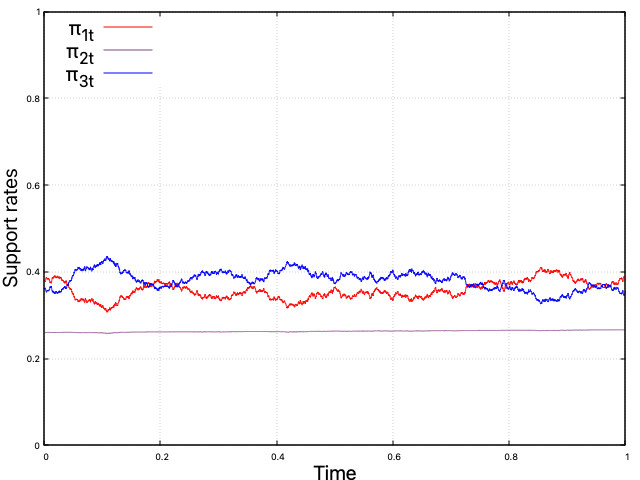}
\includegraphics[width=0.45\textwidth]{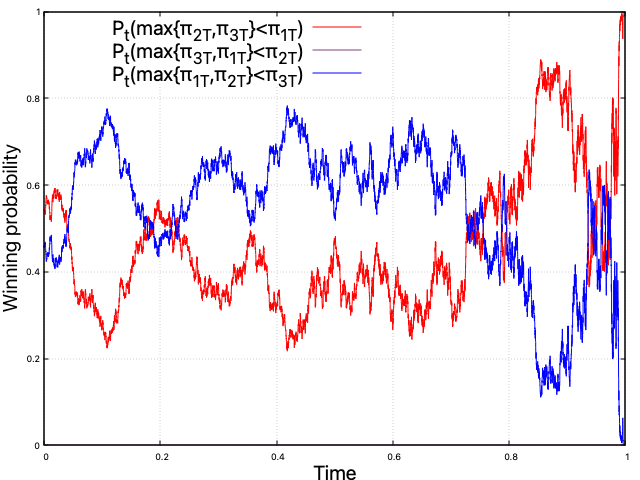} }
\centerline{
\includegraphics[width=0.45\textwidth]{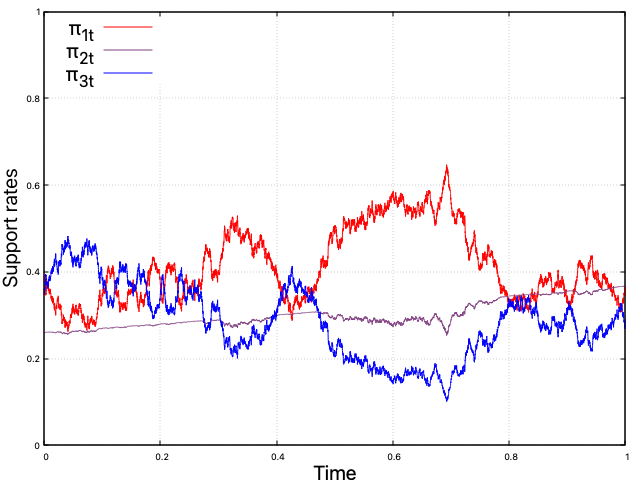}
\includegraphics[width=0.45\textwidth]{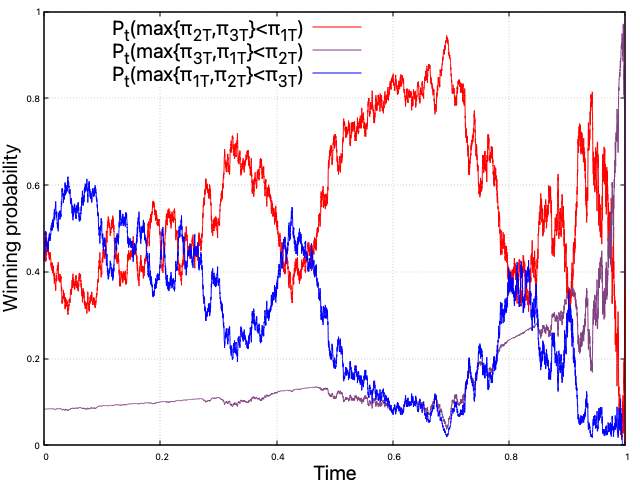}
}
\caption{\footnotesize{
\textit{Dynamical behaviours of the poll statistics $\{\pi_{it}\}$ and the 
corresponding winning probabilities}. Sample paths for the support 
rates $(\pi_{1t},\pi_{2t},\pi_{3t})$ for the three candidates are shown 
on the left panels. The corresponding winning probability processes for 
each candidate are shown on the right panels. The parameters are 
chosen as $(x_1,x_2,x_3)=(1,2,3)$ for the values of the random variable 
$X$, $(p_1,p_2,p_3)=(0.38,0.26,0.36)$ for the current support level so 
that the electorates are slightly polarised, 
and $T=1$ year for the time left to the election day. The top two panels 
correspond to the value $\sigma=0.25$ for the information flow rate. In 
this case, the probability for the second candidate to win the election is 
identically zero. For a comparison, the corresponding results for the 
choice $\sigma=1$ are plotted in the bottom two panels, in which the 
second candidate narrowly secures a victory. 
}}
\label{fig:3}
\end{figure}

In Figure~\ref{fig:2} we sketch the behaviours of the winning probabilities 
for the three candidates as functions of the current support rates $(p_1,p_2)$ 
for the first two candidates. One distinguishing feature here, as compared to 
the results for two candidate scenario, is the dependence on the information 
flow rate $\sigma$. In the two candidate case, for a given level $(p,1-p)$ of 
current support, the winning probabilities are either increasing or decreasing 
in $\sigma$. That is, if the candidate is leading the poll, then it is best not to 
reveal information, and conversely for the other candidate. In the three 
candidate case, depending on the level $(p_1,p_2,1-p_1-p_2)$ of current 
support, the winning probabilities can lack monotonicity. That is, there are 
values of $(p_1,p_2,1-p_1-p_2)$ for which a candidate will benefit from, 
say, increasing the information flow rate slightly to enhance the probability 
of winning the future election, but if it is increased too much, then this will 
result in decreasing the probability again. It follows that in a three candidate 
race, the optimal strategy of controlling information can be quite nontrivial 
for certain values of the current support rates $(p_1,p_2,1-p_1-p_2)$. 

Although we have worked out here the initial (time $t=0$) probability of a 
given candidate winning the election, it is straightforward to work out the 
corresponding conditional probabilities, such as 
${\mathbb P}(\pi_{1T}<\pi_{2T}<\pi_{3T}| \xi_t)$ and so on, that depend 
on how information has been unravelled up to time $t$. Then we are 
able to simulate not just the support rates $\{\pi_{it}\}$ but also the 
realised winning probabilities, as illustrated in Figure~\ref{fig:3}.

\section*{Representing political spectrum}

It is important to emphasise the fact that while in a two-candidate electoral 
competition the choice of the labelling numbers $\{x_0,x_1\}$ can be made 
arbitrarily, this is no longer the case if there are more than two candidates. 
The reason can be explained as follows. 

We note first that from a signal detection perspective, whether the signal 
is an unknown random variable $X$ or a (known) constant addition of an 
unknown $X$, there is no difference in the inference. This follows from the 
fact that 
the filtration generated by $\sigma X t + B_t$ is identical to that generated 
by $\sigma (X+c) t + B_t$, where $c$ is a known constant. Hence the only 
quantities of relevance in the choice of $X$ are the gaps $\omega_{ij}=
x_j-x_i$. Now in the binary case there is only one such gap $\omega=
x_1-x_0$, but scaling the gap according to $\omega\to\lambda\omega$ is 
entirely equivalent to scaling the information flow rate $\sigma\to
\lambda\sigma$. Putting it differently, the scaling 
$\omega\to\lambda\omega$ can be compensated by the scaling 
$\sigma\to\lambda^{-1}\sigma$ so that $\omega$ can be chosen arbitrarily 
by regarding $\sigma$ as the variable parameter. 

In the case of an election with three candidates, there are three gaps 
$\omega_{12}$, 
$\omega_{23}$, and $\omega_{31}$, with one constraint $\omega_{12}+
\omega_{23}+\omega_{31}=0$. Hence there are two independent scaling 
parameters, which cannot be simultaneously absorbed by scaling $\sigma$. 
It follows that the probability of a given candidate winning the election is 
dependent on the choice of the gaps $\{\omega_{ij}\}$. 
Alternatively stated, there is a natural ordering (i.e. spectrum) encoded in 
the random variable $X$ representing the candidates. In particular, the 
three candidates cannot be placed on an equal footing, for, while it is 
possible to set $\omega_{12}=\omega_{23}$, it is not possible to set  
$\omega_{12}=\omega_{23}=\omega_{31}$.  

We can take advantage of this feature of the model by observing that there 
is a long established notion of a `political spectrum' in an electoral process, 
and we can encode this information naturally in the choices of the gaps 
$\{\omega_{ij}\}$. Thus, for example, if candidate $1$ is on the left, 
candidate $2$ is moderately on the right, and candidate $3$ is further 
on the right, 
then we can let $\omega_{12}>\omega_{23}$ to capture this composition; 
and similarly for other situations. Realising this, we see that the choice of 
the gaps $\{\omega_{ij}\}$ is not up to the modeller, but it is up to the 
candidates in terms of where they place themselves in the political spectrum. 

With this in mind, we find that there are circumstances in which taking the 
political centre ground leads to a disadvantage. This follows from the 
observation that the probability for candidate $2$ to win the election is 
identically zero if $z_{12}>z_{31}>z_{23}$, while such a constraint does not 
exist for the candidate to the left or to the right. Note however that this does 
not mean that taking the centre ground is always disadvantageous --- it merely 
shows that in certain situations, having popular competitors to both the left 
and the 
right can be fatal. In particular, such a trap for candidate $2$ can be created 
among a politically polarised set of voters so that $p_2<p_1,p_3$ holds 
while at the same time $p_1\sim p_3$. 

\begin{figure}[t!]
\centerline{
\includegraphics[width=0.75\textwidth]{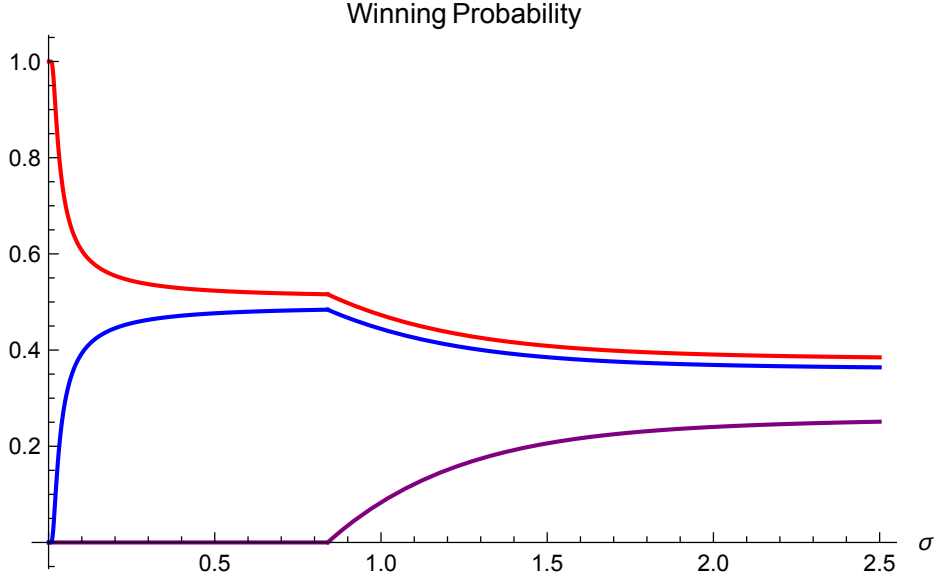}
}
\caption{\footnotesize{
\textit{Winning probabilities as functions of $\sigma$}. 
The probability of winning an election in one year time ($T=1$), as a 
function of the information flow rate $\sigma$, is 
shown for the three candidates, labelled according to $x_1=1$, 
$x_2=2$, and $x_3=3$. The current poll statistics are taken to be 
$p_1=0.38$ for the first candidate on the left (red), $p_2=0.26$ for the 
second candidate taking the centre ground (purple), and $p_3=0.36$ 
for the third candidate on the right (blue). 
}}
\label{fig:4}
\end{figure}

More specifically, a calculation shows that the condition 
$z_{12}>z_{31}>z_{23}$ can be translated into a bound on the information 
flow rate $\sigma$ as follows: 
\begin{eqnarray}
\sigma^2 < \frac{2}{\omega_{12}\omega_{23}\omega_{31}T} \, \min \left\{ 
\omega_{12} \log\frac{p_1}{p_3}-\omega_{13} \log\frac{p_1}{p_2} 
\, , \, 
\omega_{32} \log\frac{p_1}{p_3} - \omega_{13} \log\frac{p_3}{p_2} \right\} .
\label{eq:x23} 
\end{eqnarray} 
That is, provided that the inequality (\ref{eq:x23}) holds, the probability of 
candidate $2$ winning the election is identically zero. This situation is 
illustrated in Figure~\ref{fig:4}. A closer inspection 
shows that if $p_2$ is small, then the bound on $\sigma$ can be large. It 
follows that among a politically polarised electorate, the only way in which a 
candidate holding the centre ground has any chance of winning the election 
is to ensure that a lot of reliable information is revealed so as to increase the 
volatility of the poll statistics $\{\pi_{it}\}$. 

\begin{figure}[t!]
\centerline{
\includegraphics[width=0.3\textwidth]{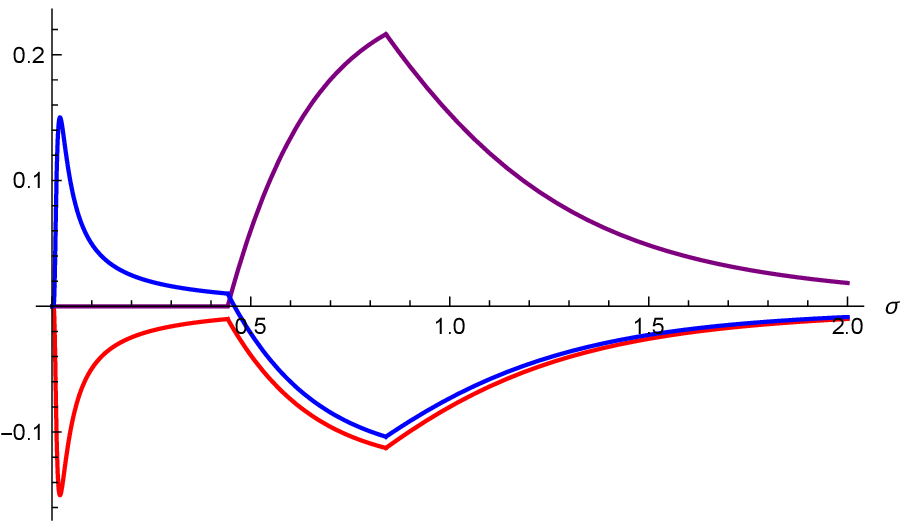}
\includegraphics[width=0.3\textwidth]{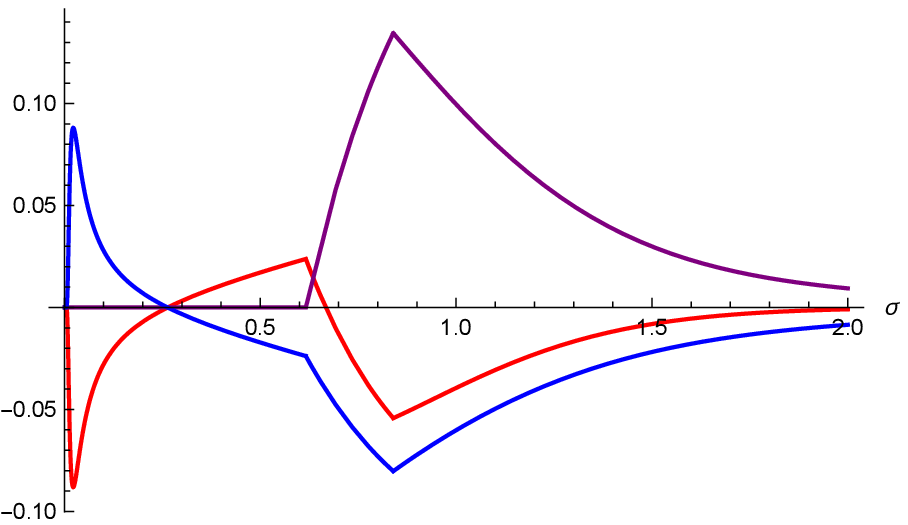}
\includegraphics[width=0.3\textwidth]{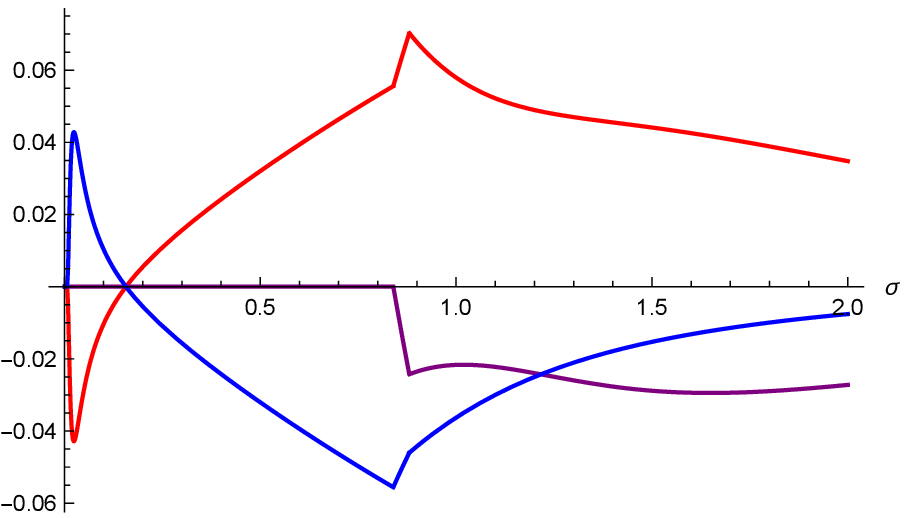}
}
\caption{\footnotesize{
\textit{Gains in winning probabilities as functions of $\sigma$}. 
If the political positioning $(x_1,x_2,x_3)=(1,2,3)$ considered in 
Figure~\ref{fig:4} is shifted, how would that affect the winning 
probabilities? Here, the difference of the resulting winning probabilities 
to the one in Figure~\ref{fig:4} is shown for three different cases: 
$(x_1,x_2,x_3)=(0.1,2,3.9)$ (left panel), $(x_1,x_2,x_3)=(1,2,3.9)$ 
(central panel), and $(x_1,x_2,x_3)=(1.5,2,3.9)$ (right panel). Other 
parameters are kept unchanged ($p_1=0.38$, $p_2=0.26$, 
$p_3=0.36$, and $T=1$). If the difference is negative, then clearly 
the shift is disadvantageous. The result shows that among a polarised 
electorate, if the candidate on the left of the political spectrum leans 
further to the left and the candidate on the right leans further to the 
right, then this is generally disadvantageous for both. However, if the 
competition is dominated by noise (small $\sigma$ values), then the 
candidate on the right can benefit by leaning further to the right. 
}}
\label{fig:5}
\end{figure}

We now ask a related question on positioning within a political spectrum. 
For this purpose we shall take the convention that if $x_j<x_k$ then 
candidate $j$ is placed politically to the left of candidate $k$. The question 
that we ask here more specifically is whether the winning probability can be 
enhanced by leaning further to the left or to the right. The answer will be 
dependent on the various parameter values, but let us consider the 
politically polarised case as shown in Figure~\ref{fig:4}. In this case, if we 
keep the value of $x_2$ unchanged but increase $x_3$ and simultaneously 
decrease $x_1$, then 
we find that the probability of winning the election for the candidate on the 
left decreases for all values of $\sigma$. However, for the candidate on the 
right, the situation is a little more complex. When the election process is 
overshadowed by noise (i.e. small $\sigma$ values), the winning probability 
of the candidate on the right can be enhanced considerably by leaning further 
to the right; whereas if the election process is not dominated by noise, then the 
winning probability decreases by leaning further to the right. Hence in this 
scenario there is no advantage for the candidate on the left to lean further 
to the left, but the candidate on the right has the advantage of turning more 
extreme, provided that the noise level is sufficiently high. If however the 
candidate misjudges the level of noise, then this strategy will backfire. Some 
examples illustrating this feature are illustrated in Figure~\ref{fig:5}, which 
shows, for example, that if the candidate on the right leans further to the 
right, while 
the candidate on the left leans slightly to the right, then there is a significant 
benefit to the candidate on the left, provided that the level of noise is not 
overwhelming.

\section*{Calibration and implementation}

As regards the practical implementation of the model, let us examine the 
model parameters that can be calibrated, and those that can be controlled. 
To this end, the current support rates $\{p_j\}$ can be fixed from today's 
poll, whereas $T$ is fixed by the date of the election. The variables 
$\{x_j\}$ that the random variable $X$ can take can then be fixed, up to 
an overall scale that can be absorbed into $\sigma$, by the relative 
positioning of the candidates within the political spectrum. A candidate, 
in particular, will have the choice for their own value of the number $x_k$, 
but will not have control over the positioning of other candidates. 

Finally, as for the value of the information flow rate $\sigma$, its value 
cannot be controlled by any individual, but its current value can be 
estimated by studying the time series for the poll statistics. This 
follows from the fact that the volatility of the support rates $\{\pi_{it}\}$ 
is given by $\sigma$. Hence a historic estimate can be used to fix 
the value of $\sigma$. Alternatively, from the odds of election betting 
it is possible to work out the implied volatility, which can be used to fix 
the value of $\sigma$. 

Having fixed all the parameters, the model can be used to interpolate 
the statistics from today to the election day. If the resulting probability 
of a given candidate winning the election appears undesirable, then 
that candidate will have a small room to manoeuvre so as to increase 
the success probability. First, the candidate can adjust their positioning 
within the political spectrum. Of course other candidates may also 
adjust their positions as a consequence of this, but if everything else 
being the same, our formula shows in which way the candidate should 
position themselves within the political spectrum. 

The other variable that a candidate can adjust is the information flow 
rate. While no individual can fix the value of $\sigma$, the result of 
(\ref{eq:4}) shows in which way an individual flow rate $\sigma_k$ 
affects the overall value of $\sigma$. In particular, $\sigma$ is a 
monotonic function of each $\sigma_k$, so increasing the value of 
any one of $\sigma_k$ will increase the overall value of $\sigma$, 
and similarly decreasing the value of $\sigma_k$ will decrease 
$\sigma$. Again, other candidates, or other information source such 
as the press, may adjust their information revelation rate as a 
result to counterbalance the impact. This, however, is just a fact 
about a democratic process -- no one candidate can control its 
outcome. Nevertheless, our framework offers an immediately 
implementable procedure for guiding the candidates to identify 
which informational strategy will increase their chances of success, 
if everything else remained the same.

\section*{Discussion}

We have examined in some detail the probability of a candidate winning 
a future election and how it is affected by control variables such as the 
level of information revelation, or noise, and the positioning of the 
candidates within the political spectrum, in the case of an electoral 
competition involving three candidates. It should be evident that a closed 
form expression for a given candidate winning a future election can be 
obtained when there are more than three candidates. For example, if 
there are four candidates, then there are 24 different ways in which the 
support rates for the candidates on the election day can be ordered, e.g., 
$\pi_{2T}<\pi_{3T}<\pi_{1T}<\pi_{4T}$ and so on. Each one of these will 
give rise to a bound on the random variable $\xi_T$ in the form of 
$\xi_T\in{\cal D}_{2314}$ for some domain ${\cal D}_{2314}$ on the real 
line. (The analogue of these domains in the case of a three-candidate 
electoral competition would be ${\cal D}_{231}=[z_{23},z_{31}]$, and so 
on.) The probability of this event being realised is therefore given by 
\begin{eqnarray}
{\mathbb P}\left(\pi_{2T}<\pi_{3T}<\pi_{1T}<\pi_{4T}\right) = 
{\mathbb E}^{{\mathbb Q}}\left[\sum_{j=1}^{4}p_{j}\exp\left( 
\sigma x_{j}\xi_{T}-\half\sigma^{2}x_{j}^{2}T\right) 
{\mathds 1}\{\xi_{T}\in{\cal D}_{2314}\}\right], 
\end{eqnarray} 
but because $\xi_T$ under ${\mathbb Q}$ is Gaussian with mean zero 
and variance $T$, this expectation can easily be worked out. By 
repeating the procedure for the five other domains ${\cal D}_{2134}$, 
${\cal D}_{3124}$, ${\cal D}_{3214}$, ${\cal D}_{1234}$, and 
${\cal D}_{1324}$, and adding the results, we obtain the probability of the 
fourth candidate winning the election. Evidently, an analogous calculation 
can be performed for each of the other candidates to work out their 
success probabilities. 

One of the nontrivial features that emerges when the number of candidates 
is greater than two is that there is a disadvantage for candidates positioning 
in the middle of the political spectrum, in a situation where the voters do not 
have strong preferences on centre grounds. An analogous property is seen 
also in the structural approach. The mathematical reason underlying this 
feature is as follows. If we label the $N$ candidates such that 
$x_1<x_2<\cdots<x_N$, then for each fixed $T$ and $\sigma$ we find that 
$\pi_{1T}(\xi_T)$, viewed as a function of $\xi_T$, is monotonically 
decreasing in $\xi_T$ without bound in the range $[0,1]$ and 
$\pi_{NT}(\xi_T)$ is monotonically increasing in $\xi_T$ without bound in 
the range $[0,1]$. However, for any $k\neq 1,N$, the function 
$\pi_{kT}(\xi_T)$, which gives the support rate for the $k$th candidate on 
the election day, is unimodal and has a maximum value at $\xi_T=\xi_k^*$ 
that is strictly less than one, where $\xi_k^*$ is the unique solution to the 
equation  
\begin{eqnarray}
\sum_{j=1}^N (x_j-x_k) p_j \exp\left( \sigma x_j \xi_k^*-\half \sigma^2 
x_j^2 T\right) = 0.
\label{eq:x25}
\end{eqnarray}

\begin{figure}[t!]
\centerline{
\includegraphics[width=0.3\textwidth]{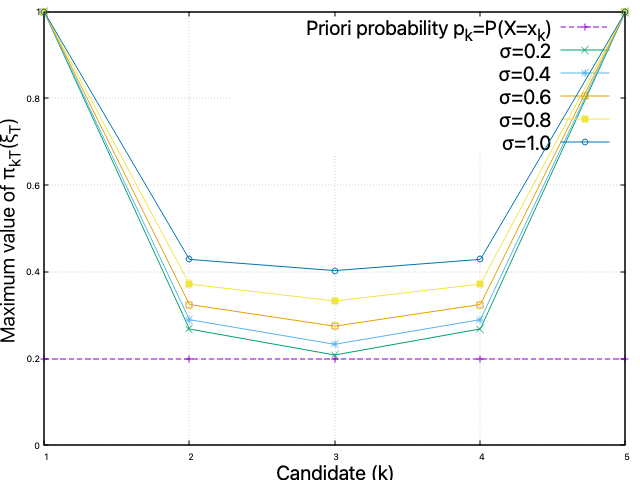}
\includegraphics[width=0.3\textwidth]{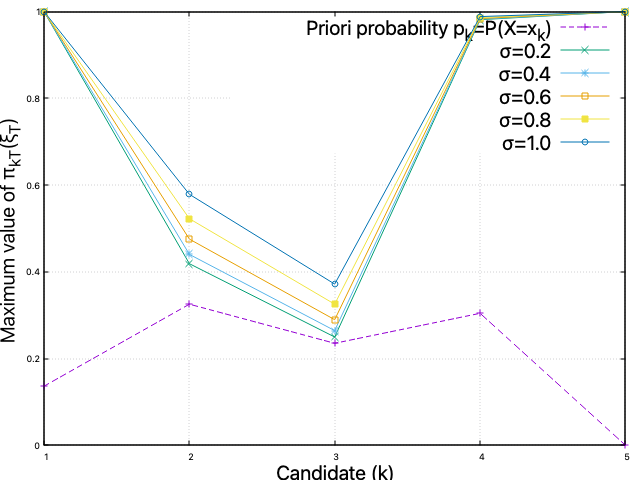}
\includegraphics[width=0.3\textwidth]{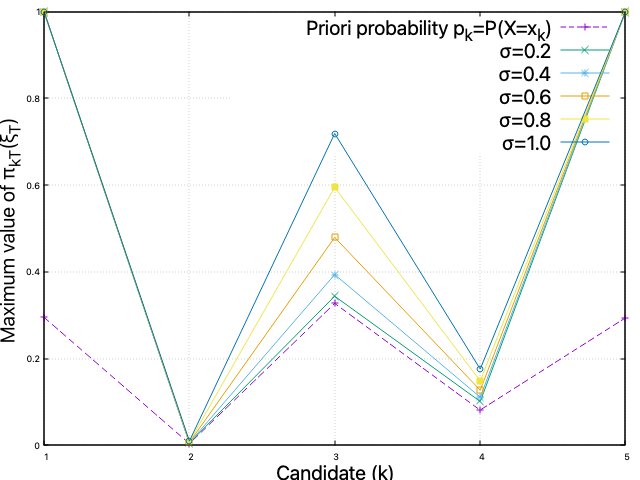}
}
\caption{\footnotesize{
\textit{Maximum attainable support rates $\pi_{kT}(\xi_k^*)$ for the five 
candidates}. 
For a range of values for the information flow rate $\sigma$, the maximum 
values of $\{\pi_{kT}\}$ on the election day 
are shown by the dots, interpolated by lines to make 
the comparison easy. The current support rates $\{p_k\}$ are given by the 
bottom values (in purple). In the left panel the five candidates are all 
assumed to have an equal support rate of 20\%, whereas they are chosen 
at random in the central and right panels. The results shows how the 
maximum attainable support rates for different candidates vary rather 
dramatically, depending on the existence of candidates having different 
political leanings and their associated current support rates. The parameters 
are chosen to be $(x_1,x_2,x_3,x_4,x_5)=(1,2,3,4,5)$ for the positioning 
of the candidates and $T=1$ year for time left to the election.  
}}
\label{fig:6}
\end{figure}

If at least one of the variables $p_k$, $\sigma$, or $T$ is large, then the 
upper bound on $\pi_{kT}$ will be close to one, so there is little concern 
for the candidate, but otherwise, the upper bound can be smaller than 
$1/N$. In the latter case, whatever information is to be circulated, the 
probability of the $k$th candidate winning a future election is identically 
zero. When there are many candidates, the threshold value $1/N$ is small, 
but if there are only three or four candidates then this effect is highly 
nontrivial and should not be ignored. Indeed, as we have seen in the case 
of a three candidate race in Figure~\ref{fig:4}, there is a wide range of values 
for the information flow rate $\sigma$ for which there is no chance for the 
second candidate to win the election. In Figure~\ref{fig:6} we show the 
maximum attainable support rates for each candidate in the case of an 
election with five candidates. The central panel in Figure~\ref{fig:6} shows 
that while the initial support rates $p_2$ and $p_4$ for the centre left and 
centre right candidates are very close, the existence of a far-right candidate 
with negligible current support level implies that the maximum attainable 
value of $\pi_{4T}$ is close to $1$, whereas the existence of a far-left 
candidate with a moderate current support level implies that  the maximum 
attainable value of $\pi_{2T}$ is considerably lower than that of $\pi_{4T}$. 

Another nontrivial feature that emerges when there are more than two 
candidates is the subtle dependence of the winning probability on the 
positioning of the candidates within the political spectrum. We have merely 
uncovered for illustration a small number of features shown in 
Figure~\ref{fig:5} and Figure~\ref{fig:6}, but a detailed sensitivity study of 
the winning probabilities on the spectrum $\{x_j\}$ is entirely feasible on 
account of the fact that we have closed-form expressions for these 
probabilities. 

We conclude by remarking how our model might be extended. Throughout 
the paper we have assumed that the information revelation rate $\sigma$ 
is constant, but in reality this is never the case. The analysis presented here 
nevertheless extends naturally to the case in which $\{\sigma_t\}$ has a 
deterministic time dependency. Specifically, in this case the candidate 
support rates take the form 
\begin{eqnarray}
\pi_{it} =\frac{p_i\exp\left( x_i \int_0^t \sigma_s \rd\xi_s-
\frac{1}{2} x_i^2 \int_0^t \sigma_s^2 \rd s\right)} {\sum_j p_j 
\exp\left( x_j \int_0^t \sigma_s \rd\xi_s-
\frac{1}{2} x_j^2 \int_0^t \sigma_s^2 \rd s\right)} .
\label{eq:x26}
\end{eqnarray}
Then the relevant random variable $\xi_T$ is replaced with $\int_0^T 
\sigma_s \rd\xi_s$, which remains Gaussian under ${\mathbb Q}$, so 
the various probabilities can still be worked out while taking into account 
a time-varying informational strategy $\{\sigma_t\}$. In this way, impacts 
of a range of time-dependent informational strategies can be 
investigated. Of course, the information process $\{\xi_t\}$ is meant to 
represent the aggregate of the various information sources, and so is  
the variable $\sigma$ (see \cite{Brody4} for how the aggregated 
information flow rate is related to that of the individual information source); 
whether it is time dependent or not. 
Thus no one candidate has the access to control the overall value of 
$\sigma$. Nevertheless, the idea is that each candidate can influence 
the value of $\sigma$, which in turn will modify the likelihoods of the 
candidates winning the election.

\vspace{0.45cm} 
\noindent 
{\bf Acknowledgements}.
The authors thank Bernhard Meister for stimulating discussion. 
DCB acknowledges support from the EPSRC 
(EP/X019926) and the John Templeton Foundation (62210). The opinions 
expressed in this publication are those of the authors and do not necessarily 
reflect the views of the John Templeton Foundation. TY is supported by 
JSPS KAKENHI (22K13965).


\begin{thebibliography}{}

\bibitem{Brody1} 
Brody,~D.C. \& Meier,~D.~M. (2022) 
Mathematical models for fake news. In 
{\em Financial Informatics: An Information-Based Approach to Asset Pricing}, 
D.~C.~Brody, \textit{et al}. (eds) (Singapore: World Scientific). 
(First appeared in 2018 in \url{https://arxiv.org/abs/1809.00964})  

\bibitem{Brody2} 
Brody,~D.C. (2019) 
Modelling election dynamics and the impact of disinformation. 
{\em Information Geometry} \textbf{2}, 209-230. 
\url{doi.org/10.1007/s41884-019-00021-2} 

\bibitem{Brody3} 
Brody,~D.~C. (2022) 
Noise, fake news, and tenacious Bayesians. 
{\em Frontiers in Psychology} \textbf{13}, 797904. 
\url{doi.org/10.3389/fpsyg.2022.797904}

\bibitem{xxx} 
Davis,~O.~A., Hinich,~M.~J. \& Ordeshook,~P.~C. (1970) 
An expository development of a mathematical model of the electoral process.  
{\em The American Political Science Review} \textbf{64}, 
426-448. 

\bibitem{Shepsle} 
Shepsle,~K.~A. (1972)  
The Strategy of Ambiguity: Uncertainty and electoral competition.  
{\em The American Political Science Review} \textbf{66}, 
555-568. 
\url{doi.org/10.2307/1957799}

\bibitem{JH} 
Harrington,~J.~E.,~Jr. (1982) 
Modelling the role of information in elections. 
{\em Mathematical and Computer Modelling} \textbf{16}, 133-145. 

\bibitem{MO} 
McKelvey,~R.~D. \& Ordeshook,~P.~C. (1985) 
Elections with limited information: A fulfilled expectations model using 
contemporaneous poll and endorsement data as information sources. 
{\em Journal of Economic Theory} \textbf{36}, 55-85. 

\bibitem{Coughlin} 
Coughlin,~P.~J. (1992) 
{\em Probabilistic Voting Theory}. 
(Cambridge: Cambridge University Press). 

\bibitem{FP} 
Feddersen,~T. \& Pesendorfer,~W. (1997) 
Voting behaviour and information aggregation in elections with private 
information. 
{\em Econometrica} \textbf{65}, 1029-1058. 

\bibitem{FM} 
Fowler,~A. \& Margolis,~M. (2014) 
The political consequences of uninformed voters. 
{\em Electoral Studies} \textbf{34}, 100-110. 

\bibitem{Lloyd} 
Rowden,~J., Lloyd,~D.~J.~B. \& Gilbert,~N. (2014) 
A model of political voting behaviours across different countries. 
{\em Physica} A\textbf{413}, 609-625. 

\bibitem{yyy} 
Coughlin,~P.~J. (2015) 
Probabilistic voting in models of electoral competition. 
In {\em Handbook of Social Choice and Voting}, J.~C.~Heckelman \& 
N.~R.~Miller (eds.) (Cheltenham: Edward Elgar Publishing Ltd.). 

\bibitem{YAKM} 
Yang,~V.~C., Abrams,~D.~M., Kernell,G. \& Motter,~A.~E. (2020) 
Why are U.S. parties so polarized? A “satisficing” dynamical model. 
{\em SIAM Review} \textbf{62}, 646-657. 
\url{doi.org/10.1137/19M1254246}

\bibitem{BHM} 
Brody,~D.~C., Hughston,~L.~P. \& Macrina,~A. (Eds.) (2022). 
{\em Financial Informatics: An Information-Based Approach to Asset 
Pricing}. (Singapore: World Scientific).

\bibitem{wonham}
Wonham,~W.~M. 1965 
Some applications of stochastic differential equations to optimal nonlinear filtering.
\textit{Journal of the Society for Industrial and Applied Mathematics} 
A\textbf{2}, 347-369.
\url{doi.org/10.1137/030202}

\bibitem{Buonaguidi} 
B. Buonaguidi,~B. 2023 
An optimal sequential procedure for determining the drift of a Brownian 
motion among three values. 
{\em Stochastic Processes and their Applications} \textbf{129}, 320-349. 
\url{doi.org/10.1016/j.spa.2023.02.001}

\bibitem{Ekstrom} 
Ekström,~E. \& Vaicenavicius,~J. 2015 
Bayesian sequential testing of the drift of a Brownian motion. 
{\em ESAIM Probab. Stat.} \textbf{19}, 626-648. 
\url{doi.org/10.1051/ps/2015012}

\bibitem{Brody4} 
Brody,~D.~C. \& Law,~Y.~T. (2015) 
Pricing of defaultable bonds with random information flow. 
{\em Applied Mathematical Finance} \textbf{22}, 399-420. 
\url{doi.org/10.1080/1350486X.2015.1050151}



\end{thebibliography}
\end{document}